\documentclass[12pt,leqno]{amsart}
\usepackage{amssymb}
\usepackage{amsmath,amssymb}
\renewcommand{\baselinestretch}{1.1}

\begin{document}
\theoremstyle{plain}
\newtheorem{MainThm}{Theorem}
\newtheorem{thm}{Theorem}[section]
\newtheorem{clry}[thm]{Corollary}
\newtheorem{prop}[thm]{Proposition}
\newtheorem{lem}[thm]{Lemma}
\newtheorem{deft}[thm]{Definition}
\newtheorem{hyp}{Assumption}
\theoremstyle{definition}
\newtheorem{rem}[thm]{Remark}
\newtheorem*{acknow}{Acknowledgments}
\numberwithin{equation}{section}
\newcommand{\eps}{{\varphi}repsilon}
\renewcommand{\d}{\partial}
\newcommand{\re}{\mathop{\rm Re} }
\newcommand{\im}{\mathop{\rm Im}}
\newcommand{\R}{\mathbb{R}}
\newcommand{\C}{\mathbf{C}}
\newcommand{\N}{\mathbf{N}}
\newcommand{\D}{C^{\infty}_0}
\newcommand{\HHH}{\mathcal{H}}
\renewcommand{\O}{\mathcal{O}}
\newcommand{\dbar}{\overline{\d}}
\newcommand{\supp}{\mathop{\rm supp}}
\newcommand{\abs}[1]{\lvert #1 \rvert}
\newcommand{\csubset}{\Subset}
\newcommand{\detg}{\lvert g \rvert}
\newcommand{\dd}{\mbox{div}\thinspace}
\newcommand{\www}{\widetilde}
\newcommand{\ggggg}{\mbox{\bf g}}
\newcommand{\ep}{\varepsilon}
\newcommand{\la}{\lambda}
\newcommand{\va}{\varphi}
\newcommand{\ppp}{\partial}
\newcommand{\ooo}{\overline}
\newcommand{\OOO}{\Omega}
\newcommand{\OOOD}{\Omega_{\delta}}
\newcommand{\wwwt}{\widetilde{t}}
\newcommand{\xxxj}{x_0^j}
\newcommand{\xxxjp}{x_0^{p(j)}}
\newcommand{\xxxjpp}{x_0^{p(j+1)}}
\newcommand{\tttj}{t_j}
\newcommand{\tttmj}{t_{j-1}}
\newcommand{\sumkj}{\sum_{k,j=1}^n}
\newcommand{\walpha}{\widetilde{\alpha}}
\newcommand{\wbeta}{\widetilde{\beta}}
\newcommand{\weight}{e^{2s\va}}
\newcommand{\fdif}{\partial_t^{\alpha}}
\newcommand{\LLL}{L_{\lambda,\mu}}
\newcommand{\LLLL}{L_{\widetilde{\la},\widetilde{\mu}}}
\newcommand{\ppdif}[2]{\frac{\partial^2 #1}{{\partial #2}^2}}
\newcommand{\uu}{\mathbf{u}}
\renewcommand{\v}{\mathbf{v}}
\newcommand{\y}{\mathbf{y}}
\newcommand{\ddd}{\mbox{div}\thinspace}
\newcommand{\rrr}{\mbox{rot}\thinspace}
\newcommand{\Y}{\mathbf{Y}}
\newcommand{\w}{\mathbf{w}}
\newcommand{\z}{\mathbf{z}}
\newcommand{\G}{\mathbf{G}}
\newcommand{\f}{\mathbf{f}}
\newcommand{\F}{\mathbf{F}}
\newcommand{\dddx}{\frac{d}{dx_0}}
\newcommand{\CC}{_{0}C^{\infty}(0,T)}
\newcommand{\HH}{_{0}H^{\alpha}(0,T)}
\newcommand{\llll}{L^{\infty}(\Omega\times (0,t_1))}
\renewcommand{\baselinestretch}{1.5}
\renewcommand{\div}{\mathrm{div}\,}  
\newcommand{\grad}{\mathrm{grad}\,}  
\newcommand{\rot}{\mathrm{rot}\,}  

\title
[Inverse coefficient problems for a transport equation $\ldots$ 
]
{Inverse coefficient problems for a transport equation by 
local Carleman estimate}


\author[P. Cannarsa, G. Floridia, F. G\"olgeleyen, 
M. Yamamoto]{P. Cannarsa$^1$, G. Floridia$^2$, F. G\"olgeleyen$^3$, 
M. Yamamoto$^4$}
\thanks{$^1$ Piermarco Cannarsa, Department of Mathematics, 
      University of Rome \lq\lq Tor Vergata'',
      00133 Rome, Italy,
email: {\tt cannarsa@mat.uniroma2.it}.\;\;\;
$^2$ Giuseppe Floridia, Department of Mathe\-matics and Applications \lq\lq R. Caccioppoli'',
      University of Naples Federico II,
80126 Naples, Italy, email: {\tt giuseppe.floridia@unina.it \& floridia.giuseppe@icloud.com}.\;\;\;
$^3$ Fikret G\"olgeleyen, Department of Mathematics, Bulent Ecevit University
Zonguldak, 67100 Turkey,
e-mail: {\tt  golgeleyen@yahoo.com}.\\
$^4$ Masahiro Yamamoto, Department of Mathematical Sciences, The University
of Tokyo, Komaba, Meguro, Tokyo 153, Japan,
Peoples' Friendship University of Russia (RUDN University), 
6 Miklukho-Maklaya St, Moscow, 117198, Russian Federation,
e-mail: {\tt myama@ms.u-tokyo.ac.jp}.}

\date{}

\maketitle

\begin{abstract} 
We consider the transport equation 
$\ppp_tu(x,t) + (H(x)\cdot \nabla u(x,t)) + p(x)u(x,t) = 0$ 
in $\OOO \times (0,T)$ where $\OOO \subset \R^n$ is a bounded domain, and
discuss two inverse problems which consist of determining a vector-valued 
function $H(x)$ or a real-valued function $p(x)$ by initial values 
and data on a subboundary of $\OOO$.
Our results are conditional stability of H\"older type in a 
subdomain $D$ provided that the outward normal component of $H(x)$ 
is positive on $\ppp D \cap \ppp\OOO$.
The proofs are based on a Carleman estimate where the weight function depends
on $H$.
\end{abstract}

\section{Introduction and main results}

Let $\OOO \subset \R^n$ be a bounded domain with smooth boundary and
$H = (h_1, ..., h_n) \in (C^1(\ooo{\OOO}))^n$. Let $(x\cdot y)$ and $|x|$ denote the scalar product of $x,y\in\R^n$ and the norm of $x\in\R^n,$ respectively.
Let $u \in C^1(\ooo{\OOO}\times [0,T])$ satisfy
$$
Pu(x,t): = \ppp_tu + (H(x)\cdot \nabla u) + p(x)u(x,t) = 0\quad 
\mbox{in $\OOO \times (0,T)$}                      \eqno{(1.1)}
$$
with 
$$
u(x,0) = a(x), \quad x \in \OOO,            \eqno{(1.2)}
$$
where $a\in C^1(\ooo{\OOO}).$

Fixing $p \in C(\ooo{\OOO})$, by 
$u(H,a) = u(H,a)(x,t)$, we denote one solution to (1.1) and (1.2).
We note that there exist many solutions and $u(H,a)$ is an arbitrarily 
chosen solution.

We discuss the following.\\
{\bf Inverse coefficient problem.}\\
We fix $p$.
Let $\Gamma \subset \ppp\OOO$ be a subboundary and let 
$a_1, ..., a_n$ be chosen suitably.  Then determine
$H(x) = (h_1(x), ..., h_n(x))$ near $\Gamma$ by 
boundary data $u(H,a_k)\vert_{\Gamma\times (0,T)}$, $k=1,2,..., n$.
\\
\vspace{0.1cm}

Our problem is an inverse coefficient problem where 
spatially varying coefficients should be determined by lateral 
boundary data and initial data.

An equation of type (1.1) is a transport equation and often appears
in mathematical physics.  For example, let $H(x,t)$ be the velocity 
field of some fluid and let $u$ be the fluid density.  Then  conservation 
of mass yields 
$$
\ppp_tu + \mbox{div}\, (uH) = 0 \quad \mbox{in $\OOO \times (0,T)$,}
$$
that is,
$$
\ppp_tu + (H\cdot \nabla u) + (\mbox{div}\, H)u = 0 \quad 
\mbox{in $\OOO \times (0,T)$},
$$
which corresponds to (1.1) if $H$ is independent of $t$.
We can interpret our inverse problem as the determination of 
a velocity field by boundary data on $\Gamma \times (0,T)$ and 
initial data.

In this paper, aiming at generous conditions for unknown $H$,
we mainly discuss  uniqueness and  stability locally near $\Gamma$.
More precisely, we are concerned with the existence of a subdomain 
$D \subset \OOO$ satisfying $\ppp D \cap \ppp\OOO \subset \Gamma$ where data
$u(H,a_k)\vert_{\Gamma\times (0,T)}$ with $k=1,..., n$ stably determine 
$H$ in $D$.

In order to state our main results, we introduce some notations.
Let $\nu = \nu(x)$ be the outward unit normal vector to $\ppp\OOO$ at $x$.
For arbitrarily fixed $x_0 \in \Gamma$ and constants $\delta_0>0, M>0$,
we define the admissible set $\mathcal{H}$ by
$$\left. \begin{array}{rl}
\HHH = &\{ H \in (C^1(\ooo{\OOO}))^n; \, 
\Vert H\Vert_{(C^1(\ooo{\OOO}))^n} \le M, \\
& (H(x_0)\cdot \nu(x_0)) > \delta_0, \quad 
\vert H(x)\vert > \delta_0 \quad \mbox{for $x\in \OOO$}\}.
\end{array}\right. 
                                                \eqno{(1.3)}
$$
Let $\cdot^T$ denote the transpose of a generic vector and set
$\nabla a(x) = (\ppp_1a, ..., \ppp_na)^T$ for any real-valued function $a(x)$.
\\

Henceforth, we always assume that $\ppp_tu(H_j,a_k) 
\in C^1(\ooo{\OOO}\times [0,T])$ for $j=1,2$ and $k=1,..., n$.

We are ready to state our first main result.
\\
{\bf Theorem 1 (local stability for the inverse coefficient problem).}\\
We assume
$$
\vert \mbox{det}\, (\nabla a_1(x), ..., \nabla a_n(x)) \vert \ne 0
\quad x \in \ooo{\OOO}                        \eqno{(1.4)}
$$
and  
$$\ppp_tu(H_j,a_k) 
\in C^1(\ooo{\OOO}\times [0,T]),\;
\Vert \ppp_t\nabla u(H_j,a_k)\Vert_{L^2(0,T;L^{\infty}(\OOO))} \le M_0, 
\; j=1,2, \, k=1,..., n  
                                           \eqno{(1.5)}$$
for some constant $M_0 > 0$.  
Then there exist constants $C>0$, $\theta \in (0,1)$ and 
a subdomain $D \subset \OOO$ satisfying $x_0 \in \ppp D \cap 
\ppp\OOO \subset \Gamma$ and such that
\begin{align*}
&\Vert H_1 - H_2\Vert_{L^2(D)}\\
\le &C\left(\sum_{k=1}^n \Vert \ppp_tu(H_1,a_k) 
- \ppp_t u(H_2,a_k)\Vert^{\theta}_{L^2(\Gamma \times (0,T))}
+ \sum_{k=1}^n \Vert \ppp_tu(H_1,a_k) 
- \ppp_t u(H_2,a_k)\Vert_{L^2(\Gamma \times (0,T))}\right)
\end{align*}
for $H_1, H_2 \in \HHH$.
Here the constants $C$ and $\theta$ depend on 
$T$, $\Gamma$, $\HHH$ and $M_0$, while 
$D$ depends on $T$, $\Gamma$ and $\HHH$, but it is independent of 
choices of $H_1, H_2$.
\\

Since $H(x) = (h_1(x), ..., h_n(x))$ possesses $n$ unknown components,
we need to take $n$ different measurements of the
 initial values satisfying (1.4) and 
observe boundary data.  We do not pose any assumptions on the size of 
$\Gamma$ and $T>0$, and the finiteness of the propagation speed 
of (1.1) forces  the domain $D$ in the stability result
to be in general smaller if $T>0$ is small, which can be seen from the proof.

Theorem 1 provides a conditional stability estimate of 
H\"older type for the inverse 
coefficient problem and asserts that we can estimate $H_1-H_2$ provided that  $H_1, H_2 \in \HHH$ and condition(1.5) is satisfied.
\\

Theorem 1 directly yields the following.
\\
{\bf Corollary (local uniqueness for the inverse coefficient problem).}\\
We assume (1.4).
Then there exists a subdomain $D \subset \OOO$ satisfying $x_0 \in \ppp D \cap 
\ppp\OOO \subset \Gamma$ such that, if
$$
u(H_1,a_k)\vert_{\Gamma \times (0,T)} 
= u(H_2,a_k)\vert_{\Gamma \times (0,T)}, \quad k=1,2,..., n
$$
with $H_1, H_2 \in \HHH$, then $H_1(x) = H_2(x)$ for all $x \in D$.
\\
\vspace{0.2cm}

We can prove similar uniqueness and stability results 
for the inverse problem of 
determining $p(x)$ in (1.1).
We fix $H \in \HHH$ arbitrarily and  denote  by $u(p,a)$ one solution 
to (1.1) and (1.2).  
\\
{\bf Theorem 2 (local stability for the inverse coefficient problem).}\\
For arbitrarily fixed constant $M>0$, we set 
$$
\mathcal{P} := \{ p\in L^{\infty}(\OOO);\, 
\Vert  p\Vert_{L^{\infty}(\OOO)} \le M\}.
$$
We assume 
$$
\vert a(x) \vert \ne 0, \quad x \in \ooo{\OOO},              \eqno{(1.6)}
$$
$\ppp_tu(p_j,a) \in C^1(\ooo{\OOO}\times [0,T])$, $j=1,2$, and 
$$
\Vert \ppp_t u(p_j,a)\Vert_{L^2(0,T;L^{\infty}(D))} \le M_0, 
\quad j=1,2\,,                                            \eqno{(1.7)}
$$
for some constant $M_0 > 0$.
Then there exist a subdomain $D \subset \OOO$ satisfying $x_0 \in \ppp D \cap 
\ppp\OOO \subset \Gamma$ and constants $C>0$ and $\theta \in (0,1)$, 
depending on $T$, $\Gamma$, $\HHH$ and $M_0$, such that 
\begin{align*}
& \Vert p_1 - p_2\Vert_{L^2(D)}\\
\le &C(\Vert \ppp_tu(p_1,a) - \ppp_t u(p_2,a)\Vert^{\theta}
_{L^2(\Gamma \times (0,T))}
+ \Vert \ppp_tu(p_1,a) - \ppp_t u(p_2,a)\Vert_{L^2(\Gamma \times (0,T))})
\end{align*}
for $p_1, p_2 \in \mathcal{P}$ satisfying (1.7).
Here $D$ depends on $T$, $\Gamma$ and $\HHH$, but it is independent of 
the choice of $p_1, p_2$.
\\
\vspace{0.1cm}

In Theorems 1 and 2, we do not know boundary values on $\ppp D \setminus 
\Gamma$, and Lemma 2 in Section 2 proves that $(H\cdot\nu) \le 0$ on 
$\ppp D \setminus \Gamma$, where $\nu$ is the outward unit normal vector 
to $D$.  Therefore, in terms of data on the inflow and outflow subboundaries,
we note that we cannot uniquely determine
the solutions $u(H_j,a_k)$, $u(p_j,a)$ themselves on the whole domain
$D\times (0,T)$ although the coefficients are uniquely determined in $D$.
In fact, for the uniqueness of solution $u$ in 
$D \times (0,T)$, data of $u$ on $\{ x\in \ppp D; \, 
H(x)\cdot \nu(x) \le 0\}$ are indispensable.
We give an illustrating simple example:
\\
{\bf Example.}
$$
\left\{ \begin{array}{rl}
& \ppp_tu(x,t) - \ppp_xu(x,t) = 0, \quad 0<x,t<1, \\
& u(x,0) = 0, \quad 0 < x < 1, \\
& u(0,t) = 0, \quad 0 < t < 1.
\end{array}\right.
$$
We can regard $u$ as the difference of two solutions 
of (1.1) with $p=0$ satisfying (1.2) with same $a$.
Let $g \not\equiv 0$ be a function of class $C^1(\R)$ satisfying $g(\eta) = 0$ for $\eta \le 1$.
Then we can easily see that $u(x,t) = g(x+t)$ is a solution to the above
initial-boundary value problem in $(0,1)^2$. This example 
indicates that no uniqueness result can be expected for the solution $u$ in the whole domain.
Here $H(0) = -1$ and 
$\nu(0) = -1$ at $x=0$ where boundary data are given, and we note 
that $(H(0) \cdot \nu(0)) = 1>0$, but we are 
given no boundary data on the subboundary where $(H(1)\cdot \nu(1)) < 0$.

In particular, if we know data of $u$ on $(\ppp D\setminus \Gamma) 
\times (0,T)$, then we can prove that in the inverse problem,
both coefficients and solution are uniquely determined.

Since (1.1) is a  first order equation, for the above inverse problem
one could apply the method of characteristics as was done, for example, by Belinskij \cite{Be} and Romanov 
\cite[Chapter 5]{Ro}.
However, in this paper we follow a different approach because
\begin{enumerate}
\item
the method of characteristics does not directly provide stability estimates 
in $L^2$-spaces  (like Theorems 1 and 2) for  inverse problems;
\item
for inverse coefficient problems, a more comprehensive method is available,
which is applicable not only to hyperbolic, parabolic, Schr{\"o}dinger, 
and other equations in mathematical physics, but also to the transport 
equation.
\end{enumerate}

Such a comprehensive methodology for inverse coefficient problems
was created by Bukhgeim and Klibanov \cite{BK}.  Since \cite{BK}, 
there have been many 
publications for second-order equations.  We can refer for example to 
Beilina and Klibanov \cite{BeiK}, Bellassoued and Yamamoto \cite{BY},
Imanuvilov and Yamamoto \cite{IY1} - \cite{IY4}, 
Klibanov \cite{Kli1}, \cite{Kli2}, Klibanov and Timonov \cite{KliT}, 
Klibanov and Yamamoto \cite{KliYa}, Yamamoto \cite{Ya}.
Here we do not intend to give any substantial list of 
references and the reader 
can consult the references contained in these works.
The methodology by Bukhgeim and Klibanov \cite{BK}
is based on Carleman estimates which 
are $L^2$-weighted estimates for solutions to partial differential equations.
The choice of the weight functon is crucial for each Carleman estimate.
As for Carleman estimates themselves, see H\"ormander \cite{Ho}, Isakov 
\cite{Is} for instance.

Especially for transport equations,  the approach 
based on Carleman estimates was adopted, for instance, in Cannarsa, Floridia and Yamamoto \cite{CFY}, Gaitan and Ouzzane \cite{GO}, 
and G\"olgeleyen and Yamamoto \cite{GY}.  The works \cite{GO} and 
\cite{GY} discuss  global stability in $x$ provided that 
boundary values are given on the whole lateral boundary 
$\ppp\OOO \times (0,T)$.
Our approach is similar to  G\"olgeleyen and Yamamoto \cite{GY}
but, in this paper,
by concentrating on  local stability we can allow the unknown principal part to satisfy the generous condition 
$(H(x_0)\cdot \nu(x_0)) > 0$.
Reference \cite{CFY} considers a transport equation 
$\ppp_tu + (H(t)\cdot \nabla u) = 0$ with time-dependent $H(t)$ and
establishes a Carleman estimate to prove an  observability-type inequality for the $L^2$-norm of the initial value in terms of 
boundary data.
Moreover, Klibanov and Pamyatnykh \cite{Klipa1}, \cite{Klipa2} and 
Machida and Yamamoto \cite{MY} discuss Carleman estimates and
inverse problems for a transport equation including an integral
term with respect to a velocity field.
 
In our method, the choice of the weight function in our Carleman estimate 
is essential for the application to inverse problems.
The weight function of our Carleman estimate (Proposition 1 in Section 2)
is linear in $t$ and similar to \cite{GY}, but different from \cite{GO}, 
\cite{Klipa1}, \cite{Klipa2} where the weight is quadratic in $t$.
Thanks to our choice of the weight function, 
we need not extend  the solution $u$ to $(-T,0)$. 
Moreover, as is already explained, in our Carleman estimate,
the sufficient condition $(H(x_0) \cdot \nu(x_0)) > 0$ in (1.3) seems
more related to a physical interpretation than the one in \cite{GO}.
Furthermore, for  second order hyperbolic equations one can take
even or odd  extensions of the solution $u(x,t)$  to 
establish Carleman estimates in $\OOO\times (-T,T)$.
We notice that a Carleman estimate in $(-T,T)$ is usually easier to be 
established than on $\OOO \times (0,T)$.  
For the transport equation such extensions cause a difficulty because 
the signs of the principal 
terms are not the same in $(0,T)$ and $(-T,0)$.  We note that 
by starting with (1.1) in $t \in (-T,T)$ with (1.2) and not extending in $t$,
the inverse problem is easier but this is not an inverse problem for
an initial-boundary value problem, and so it is unnatural to assume 
that the choice of $a$, $a_k$ is given a priori.   
 
To sum up, compared with the existing related papers \cite{Be}, \cite{GO},
\cite{GY}, \cite{Ro} on inverse problems for the transport equation (1.1), 
we list our achievements:
\begin{itemize}
\item
inverse problems for initial-boundary value problems;
\item
local stability in $x$: this means that we do not assume boundary data
on the whole lateral boundary $\ppp\OOO \times (0,T)$ and 
we are concerned with estimation in a small subdomain;
\item
flexible choice of the weight function for the Carleman estimate which 
essentially depends on $H(x);$
\item
estimates in $L^2$-norms unlike \cite{Be} and \cite{Ro}.
\end{itemize}
 
This paper is composed of five sections.  In Section 2, we establish 
a key Carleman estimate (Proposition 1) and 
a lemma characterizing the subdomain $D$.
Sections 3 and 4 are devoted to the proofs of Theorems 1 and 2, respectively,
and Section 5 is devoted to remarks about other possible formulations of 
inverse problems.
\section{Key Carleman estimate and a key lemma}

Let $D \subset \OOO$ be a bounded domain and let $Q = D \times (0,T)$, and
$\nu = \nu(x)$ be the outward unit normal vector to $\ppp D$ at $x$.
Without loss of generality, we assume  
$$
0 \in \ooo{D}.
$$
We set 
$$
\psi_H(x) = \sum_{j=1}^n x_jh_j(x), \quad x\in D,  \eqno{(2.1)}
$$
$$
\mu_H := \min_{x\in \ooo D} (H(x)\cdot \nabla\psi_H(x)),
                                                       \eqno{(2.2)}
$$
and
\begin{align*}
&\Gamma_+ = \{ x\in \ppp D;\thinspace (H(x) \cdot \nu(x)) \ge 0\}, \\
&\Gamma_- = \{ x\in \ppp D;\thinspace (H(x) \cdot \nu(x)) \le 0\}.
\end{align*}

We can prove  
\\
{\bf Lemma 1.}\\
Let $H \in \HHH$ be arbitrary.  If diam $D < \frac{\delta_0^2}{2M^2}$, 
then
$$
\mu := \inf_{H\in \HHH} \mu_H \ge \frac{\delta_0^2}{2}.
$$
\\
{\bf Proof.}\\
Let $H = (h_1, ..., h_n) \in \HHH$ be arbitrary.
By the Cauchy-Schwarz inequality, we have
\begin{eqnarray*}
&& (H(x)\cdot \nabla\psi_H(x)) 
= \sum_{k=1}^n h_k(x)^2 + \sum_{k=1}^nh_k(x)
\sum_{j=1}^n x_j\ppp_kh_j(x)\\
\ge &&\min_{x\in\ooo D} \vert H(x)\vert^2 
- \left(\sum_{k=1}^n h_k(x)^2\right)^{\frac{1}{2}}
\left( \sum_{k=1}^n \left\vert 
\sum_{j=1}^n x_j\ppp_{k}h_j(x)\right\vert^2\right)^{\frac{1}{2}}\\
\ge && \delta_0^2 - \Vert H\Vert_{(L^{\infty}(D))^n}
\left(\sum_{k=1}^n \left( \sum_{j=1}^n \vert x_j\vert^2\right)
\left( \sum_{j=1}^n \vert \ppp_kh_j(x)\vert^2\right)
\right)^{\frac{1}{2}}\\
\ge&& \delta_0^2 -  \max_{x\in\ooo D} \vert x\vert
\Vert H\Vert_{(L^{\infty}(D))^n}
\Vert \nabla H\Vert_{(L^{\infty}(D))^{n\times n}}
\ge \delta_0^2 - (\mbox{diam}\, D)M^2.
\end{eqnarray*}
Therefore $\inf_{H\in \HHH} \mu_H \ge \frac{\delta_0^2}{2}$ if
diam $D < \frac{\delta_0^2}{2M^2}$.
Thus the proof of Lemma 1 is complete.
\\

We further set 
$$
\va_H(x,t) = -\beta t + \psi_H(x), \quad (x,t) \in Q, \qquad
B_H(x) = (H(x)\cdot \nabla \psi_H(x)) - \beta,
$$
where $\beta>0$ satisfies $0<\beta < \mu$.
Then (2.2) implies
$$
B_H(x) \ge \mu-\beta > 0, \quad x\in \ooo{D}.   \eqno{(2.3)}
$$
\\

We can state the key Carleman estimate.
\\
{\bf Proposition 1 (Carleman estimate).}\\
Let $\ppp D$ be piecewise smooth.
There exist constants $s_0 > 0$ and
$C>0$ such that 
$$
s\int_D \vert u(x,0)\vert^2 e^{2s\va_H(x,0)} dx 
+ s^2\int_Q \vert u\vert^2 e^{2s\va_H(x,t)} dxdt
+ Ce^{-Cs}\int_{\Gamma_- \times (0,T)}
\vert (H\cdot\nu)\vert \vert u\vert^2 dSdt
$$
$$
\le C\int_Q \vert Pu\vert^2 e^{2s\va_H(x,t)} dxdt
+ Ce^{Cs}\int_{\Gamma_+\times (0,T)} \vert u\vert^2 dS dt
+ Cs\int_D \vert u(x,T)\vert^2 e^{2s\va_H(x,T)} dx \eqno{(2.4)}
$$
for all $s > s_0$, all $H\in \HHH,$ $p \in \mathcal{P}$ 
and $u \in H^1(Q).$
\\

We note that $C>0$ and $s_0>0$ are independent of choices of 
$H \in \HHH$ and $p \in \mathcal{P}$.
\\
{\bf Proof of Proposition 1.}\\
It suffices to assume that $p \equiv 0$.  Indeed let the proposition be proved 
for $p\equiv 0$.  Then, applying (2.4) to 
$\ppp_tu + (H\cdot \nabla u) = Pu - pu$, we have
\begin{align*}
& s\int_D \vert u(x,0)\vert^2 e^{2s\va_H(x,0)} dx 
+ s^2\int_Q \vert u\vert^2 e^{2s\va_H(x,t)} dxdt
+ Ce^{-Cs}\int_{\Gamma_- \times (0,T)}
\vert (H\cdot\nu)\vert \vert u\vert^2 dSdt
\\
\le &C\int_Q \vert Pu\vert^2 e^{2s\va_H(x,t)} dxdt
+ C\int_Q \vert pu\vert^2 e^{2s\va_H(x,t)} dxdt 
+ Ce^{Cs}\int_{\Gamma_+\times (0,T)} \vert u\vert^2 dS dt\\
+& Cs\int_D \vert u(x,T)\vert^2e^{2s\va_H(x,T)} dx.
\end{align*}
Since $\Vert p\Vert_{L^{\infty}(\OOO)} \le M$ by $p \in \mathcal{P}$,
we can absorb the second term on the right-hand side into the 
left-hand side by choosing $s>0$ large, so that (2.4) is seen to hold 
for general $p\in \mathcal{P}$.
\\
Let 
$$
w = e^{s\va_H}u, \quad Lw = e^{s\va_H}P(e^{-s\va_H}w).
$$
Then
$$
Lw = \ppp_tw + (H(x) \cdot \nabla w) - sB_Hw
\quad \mbox{in $Q$}.
$$
Therefore, in terms of (2.3), we have  
\begin{align*}
& \int^T_0 \int_D \vert Pu\vert^2 e^{2s\va_H} dxdt
= \int^T_0 \int_D \vert Lw\vert^2 dxdt\\
\ge & - 2s\int_Q B_Hw(\ppp_tw + (H(x)\cdot \nabla w)) dxdt
+ s^2\int_Q (\mu-\beta)^2\vert w\vert^2 dxdt\\
=-& s\int_{Q} B_H\ppp_t(\vert w\vert^2) dxdt 
- s\int_Q B_H\sum_{j=1}^n h_j\ppp_j(\vert w\vert^2) dxdt\\
+& s^2(\mu-\beta)^2\int_Q \vert w\vert^2 dxdt\\
=& s\int_D B_H(x) \left[ \vert w(x,t)\vert^2\right]^{t=0}_{t=T} dx
+ s\int_Q \mbox{div} (B_HH) \vert w\vert^2 dxdt \\
-& s\int^T_0 \int_{\ppp D} B_H(x)(H\cdot\nu) \vert w\vert^2 dSdt
+ s^2(\mu-\beta)^2\int_Q \vert w\vert^2 dxdt\\
\ge & s\int_D (\mu-\beta) \vert w(x,0)\vert^2 dx
- Cs\int_D \vert w(x,T)\vert^2 dx
+ s\int_Q \mbox{div} (B_HH) \vert w\vert^2 dxdt \\
-& s\int_{\Gamma_+\times (0,T)} B_H(x)(H\cdot\nu) \vert w\vert^2 dSdt
- s\int_{\Gamma_- \times (0,T)} B_H(x)(H\cdot\nu)\vert w\vert^2 dSdt\\
+ &s^2(\mu-\beta)^2\int_Q \vert w\vert^2 dxdt.
\end{align*}
Hence, using $\displaystyle \sup_{H\in \HHH} \Vert \va_H\Vert_{C(\ooo{D})} 
< \infty$ and
taking some constant $C(M)>0$, we obtain 
\begin{align*}
& \int_Q \vert Pu\vert^2 e^{2s\va_H} dxdt
+ Cs\int_{\Gamma_+\times (0,T)} \vert u\vert^2 e^{2s\va_H} dS dt
+ Cs\int_D \vert u(x,T)\vert^2 e^{2s\va_H(x,T)} dx\\
\ge & s\int_D (\mu-\beta) \vert u(x,0)\vert^2 e^{2s\va_H(x,0)}dx
+ s^2\left( (\mu-\beta)^2 - \frac{C(M)}{s}\right) \int_Q \vert w\vert^2 dxdt\\
+& Cs\int_{\Gamma_-\times (0,T)} \vert B_H(x)(H\cdot\nu)\vert 
\vert u\vert^2 e^{2s\va_H} dS dt.
\end{align*}
Here $C(M)>0$ is dependent on $M$, but independent of choices of 
$p \in \mathcal{P}$ and $H \in \HHH$, and 
we can choose $s_0 > 0$ and 
$C>0$ independently of the choices of $p$ and $H,$ such that (2.4) holds.
Thus the proof of Proposition 1 is completed.
\\

We conclude this section with a lemma characterizing the subdomain $D$
which can be arbitralily small as we wish.\\

{\bf Lemma 2.}
There exists a constant $r>0$ such that for $0 < \ep < r$, 
we can choose a subdomain $D = D(\ep) \subset \OOO$ 
such that $\ppp D$ is piecewise smooth, diam $D < \ep$ and
$$
\mbox{$\ppp D \cap \ppp\OOO \supset \Gamma \cap \{x \in \R^n;\,
\vert x-\www{x}\vert < \rho_0\}$}
$$
$$
\mbox{with some $\www{x} \in \R^n$ and $\rho_0>0$, 
$x_0\in \ppp D \cap \ppp\OOO$, and}
$$ 
$$
(H(x)\cdot \nu(x)) \le 0 \quad \mbox{for all $x\in \ppp D \setminus 
\Gamma$ and all $H \in \HHH$}.              \eqno{(2.5)}       
$$
\\
{\bf Proof.\\ First Step.} It suffice to prove that there exist $r>0$ and 
$D$ satisfying (2.5) and diam $D < r$.
Then by the argument below, we see that for small $\ep > 0$ satisfying
$0 < \ep < r$, there exists a domain $D(\ep)$ satisfying (2.5).

We set $x = (x', x_n) \in \R^{n-1}\times \R$,
$\nabla' \ell(x') = (\ppp_1\ell, ..., \ppp_{n-1}\ell)^T$,
$B_{\rho} = \{ x' \in \R^{n-1}; \, \vert x'\vert < \rho\}$ with 
$\rho>0$.  For a domain $E \subset \R^n$ and subboundary 
$L \subset \ppp E$, by $\nu_L(x)$ we denote the outward unit normal vector to 
$E$ at $x\in L$. 

Since we assume that $x_0=0$, from (1.3) it follows that
$$
(H(0) \cdot \nu_{\Gamma}(0)) > \delta_0 \quad \mbox{for all $H \in \HHH$}.
$$
Since $\partial\Omega$ is smooth,
there exist a constant $\rho_0>0$ and a smooth function $\ell$ 
defined in $B_{\rho_0}$ such that $\ell(0) = 0$ and 
$$
\gamma_0 = \{(x',x_n); \,x'\in B_{\rho_0}, \quad x_n=\ell(x')\}
\subset \ppp\OOO.
$$
Without loss of generality, we further assume that $\OOO$ is located below 
$x_n = \ell(x')$ locally near $0$, that is, if $\vert x'\vert$ is small
and $x_n - \ell(x')$ is positive and small, then $(x',x_n) \not\in 
\overline{\OOO}$.
Then we can represent 
$$
\nu_{\gamma}(x',x_n) = \frac{(-\nabla' \ell(x'), 1)^T}
{\sqrt{1+\vert\nabla'\ell(x')\vert^2}}.
$$
Then from the interval $(0, \min\{ \rho_0, 1\})$,
we can choose a constant $\rho_1>0$ sufficiently small 
such that 
$$
(H(x) \cdot \nu_{\Gamma}(x)) > \frac{\delta_0}{2} \eqno{(2.6)}
$$
$$
\quad \mbox{for all 
$x \in \gamma_1:= \{ (x',x_n) \in \R^n;\, x_n = \ell(x'), \,
\vert x'\vert <\rho_1\}$ and all $H \in \HHH$}.
$$
We note that $0\in \gamma_1$ and $\gamma_1 \subset \gamma_0$, and 
$\rho_1$ is independent of choices $H \in \HHH$.
\\
{\bf Second Step: Proof of (2.6).}
As is easily verified, we have 
$$
\sup_{x\in \gamma_1} \vert x\vert \le
\rho_1(1+\Vert\ell \Vert^2_{C^1(\ppp\OOO)})^{\frac{1}{2}}
=: C(M)\rho_1.
$$
Henceforth $C=C(M)>0$ denotes generic constants which are independent of 
choices of $H \in \HHH$ and $x\in \R^n$ near $0$. 
Hence
\begin{align*}
& (H(x) \cdot \nu_{\Gamma}(x)) 
= (H(0)\cdot \nu_{\Gamma}(0)) + ((H(x)-H(0)) \cdot \nu_{\Gamma}(x))
+ (H(0)\cdot (\nu_{\Gamma}(x) - \nu_{\Gamma}(0))\\
\ge &(H(0)\cdot \nu_{\Gamma}(0)) - 2M\vert x\vert\sup_{x\in \gamma_1}
\vert \nu(x)_{\Gamma}\vert 
- M\sup_{x\in \gamma_1} \vert \nu_{\Gamma}(x) - \nu_{\Gamma}(0)\vert \\
> & \delta_0 - C(M)\rho_1.
\end{align*}
Therefore if $\rho_1 < \frac{\delta_0}{2C(M)}$, then (2.6) is verified.
\\

Next we choose $\www{\ell} \in C^2(\ooo{B_{\rho_1}})$ such that 
$\www{\ell}(x')  > 0$ for $\vert x'\vert < \rho_1$ and 
$\www{\ell}(x') = 0$ for $\vert x'\vert = \rho_1$ and 
$$
\Vert \www{\ell} \Vert_{C^2(\ooo{B_{\rho_1}})} \le \frac{\delta_0}{M}, \quad 
\ppp_i\ppp_j\www{\ell}(x') = 0, \quad 1\le i,j \le n-1,\, \vert x'\vert = \rho_1.
$$
We set 
$$
D = \{ (x',x_n)\in \R^n;\, 
\ell(x') - \www{\ell}(x') < x_n < \ell(x'), \quad \vert x'\vert <\rho_1\}
$$
and 
$$
\gamma_2 = \{ (x',x_n)\in \R^n;\, 
x_n = \ell(x') - \www{\ell}(x'), \quad \vert x'\vert < \rho_1\}.
$$

Then we see that $D \setminus \ooo{\gamma_1} \subset \OOO$
and $\ppp D = \gamma_1 \cup \gamma_2$.  

Now we prove
$$
(H(x) \cdot \nu_{\gamma_2}(x)) \le -\frac{\delta_0}{4} \quad \mbox{for all 
$x \in \gamma_2$ and all $H \in \HHH$}.    \eqno{(2.7)}
$$
We note that $(H(x)\cdot \nu_{\gamma_k}(x))$, $k=1,2$ may be discontinous 
on $\ooo{\gamma_1
\cap \gamma_2}$.
\\
{\bf Third Step: Proof of (2.7).}
The outward normal vector $\nu_{\gamma_2}(x)$ to $D$ at 
$x \in \gamma_2$ is given by 
$$
\nu_{\gamma_2}(x) = \frac{(\nabla' (\ell-\www{\ell} )(x'), -1)^T}
{\sqrt{1+\vert\nabla'(\ell-\www{\ell} )(x')\vert^2}}
$$
if $\vert x'\vert < \rho_1$.  Then the definition of $\gamma_1$ and 
$\gamma_2$ yields 
$$
\nu_{\gamma_2}(x', x_n) = - \nu_{\gamma_1}(x', x_n) \quad 
\mbox{for $\vert x'\vert = \rho_1$},          \eqno{(2.8)}
$$      
and $D$ is the connected part between the two
subboundaries $\gamma_1$, $\gamma_2$.
For an arbitrarily fixed $y = (y', (\ell-\www{\ell} )(y')) \in \gamma_2$ 
with $\vert y'\vert = \rho_1$, we see from (2.6) and (2.8) that 
$$
(H(y) \cdot \nu_{\gamma_2}(y)) \le -\frac{\delta_0}{2} \quad 
\mbox{for all $H \in \HHH$}.                   \eqno{(2.9)}
$$
For each $x, y \in \gamma_2$, we have
\begin{align*}
& \vert \nu_{\gamma_2}(x) - \nu_{\gamma_2}(y)\vert 
= \biggl\vert  \left( \frac{1}{\sqrt{1+ \vert \nabla'(\ell-\www{\ell})(x')
\vert^2}}
- \frac{1}{\sqrt{1+ \vert \nabla'(\ell-\www{\ell} )(y')\vert^2}} \right)
(\nabla'(\ell-\www{\ell} )(x'), -1)^T\\
+ &\frac{1}{\sqrt{1+ \vert \nabla'(\ell-\www{\ell} )(y')\vert^2}}
(\nabla'(\ell-\www{\ell} )(x') - \nabla'(\ell-\www{\ell} )(y'), 0)^T 
\biggr\vert\\
\le& C\frac{\bigg\vert\sqrt{1+ \vert \nabla'(\ell-\www{\ell} )(x')\vert^2}
- \sqrt{1+ \vert \nabla'(\ell-\www{\ell} )(y')\vert^2} \bigg\vert}
{\sqrt{1+ \vert \nabla'(\ell-\www{\ell} )(x')\vert^2}
\sqrt{1+ \vert \nabla'(\ell-\www{\ell} )(y')\vert^2}}\\
+ & C\vert \nabla'(\ell-\www{\ell} )(x') - \nabla'(\ell-\www{\ell} )(y')\vert\\
\le& C\left\vert \vert \nabla'(\ell-\www{\ell} )(x') \vert^2
- \vert \nabla'(\ell-\www{\ell} )(y')\vert^2 \right\vert
+ C\vert \nabla'(\ell-\www{\ell} )(x') - \nabla'(\ell-\www{\ell} )(y')\vert
\end{align*}
$$
 < C\rho_1.                 \eqno{(2.10)}
$$
At the last inequality, by the mean value theorem, we argued 
\begin{align*}
& \vert \ppp_i(\ell-\www{\ell})(x') - \ppp_i(\ell-\www{\ell})(y')\vert
\le \Vert \nabla'\ppp_i(\ell-\www{\ell})\Vert_{C(\ooo{B_{\rho_1}})}
\vert x'-y'\vert\\
\le& \Vert \ell-\www{\ell}\Vert_{C^2(\ooo{B_{\rho_1}})}
(\vert x'\vert + \vert y'\vert)
\le 2C\rho_1, \quad 1\le i \le n-1
\end{align*}
for $x', y' \in B_{\rho_1}$.  
Therefore (2.9) and (2.10) imply 
\begin{align*}
& (H(x)\cdot \nu_{\gamma_2}(x))
= (H(y)\cdot \nu_{\gamma_2}(y))
+ (H(y) \cdot (\nu_{\gamma_2}(x) - \nu_{\gamma_2}(y))
+ ((H(x) - H(y)) \cdot \nu_{\gamma_2}(x))\\
< & -\frac{\delta_0}{2} 
+ \vert (H(y) \cdot (\nu_{\gamma_2}(x) - \nu_{\gamma_2}(y))\vert
+ \vert ((H(x) - H(y)) \cdot \nu_{\gamma_2}(x))\vert \\
\le & -\frac{\delta_0}{2} + CM\rho_1 + CM\rho_1.
\end{align*}
We further choose $\rho_1>0$ smaller such that 
$2CM\rho_1 < \frac{\delta_0}{4}$.
Then (2.7) is verified.
\\

We can readily see that diam $D \le \sqrt{2}\rho_1$, and setting
$r \le \sqrt{2}\rho_1$, we complete the proof of Lemma 2.

\section{Proof of Theorem 1.}

The proof is based on a key idea which is similar to 
G\"olgeleyen and Yamamoto \cite{GY}, but for our stability local in 
$x$ we need non-trivial modifications.
The proof is divided into two steps.

{\bf First Step: Cut-off and separation of the range set of the weight
function.}\\
Since $u(H_j,a_k)$ does not vanish at $t=T$, we cannot directly apply 
Proposition 1.  Therefore we need to introduce a cut-off function $\chi(t)$
vanishing near $t=T$ to consider $\chi(u(H_1,a_k) - u(H_2,a_k))$.
However the transport equations satisfied by these functions which 
are cut off by $\chi$, contain some extra terms with $\frac{d\chi}{dt}$ 
and for keeping these terms 
as minor terms, we have to verify some separation of the range set of 
the weight function $\va_{H_1}(x,t)$.

In (1.3), for simplicity, we can assume that $x_0=0$.  Let $H, H_1, H_2 \in \HHH$ be arbitrary.
We choose $\beta > 0$ satisfying 
$$
0 < \beta < \frac{\delta_0^2}{2}.
$$
We choose $\ep>0$ satisfying 
$$
0 < \ep < \min \left\{ \frac{\delta_0^2}{2M^2}, 1, \frac{\beta T}{4M}, r
\right\},
                                                      \eqno{(3.1)}
$$
where $r>0$ is the constant whose existence is guaranteed by Lemma 2.
By Lemma 2, we can construct a small domain $D$ such that 
$\ppp D \cap \ppp\OOO$ is a non-empty 
subset of the interior of
$\Gamma$, $0\in\ppp D \cap \ppp\OOO,$ diam $D < \ep$ and 
$(H(x)\cdot \nu(x)) \le 0$ for all $x \in \ppp D \setminus \Gamma$ and
all $H \in \HHH$.
Then, by (3.1) Lemma 1 implies that 
$$
\mu := \inf_{H\in \HHH} \min_{x\in \ooo{D}} (H(x)\cdot \nabla\psi_H(x))
\ge \frac{\delta_0^2}{2}.
$$
Therefore Carleman estimate (2.4) holds in $Q=D \times (0,T)$ uniformly 
with respect to $H \in \HHH$.  

Next we estimate $\displaystyle\max_{y\in \ooo{D}} \psi_H(y)
- \min_{y\in \ooo{D}} \psi_H(y)$.  We can choose $\xi = \xi(H)$,
$\eta = \eta(H) \in D$ such that 
$$
\max_{y\in \ooo{D}} \psi_H(y) - \min_{y\in \ooo{D}} \psi_H(y)
= \psi_{H}(\xi) - \psi_H(\eta).
$$
Hence, since $\Vert H\Vert_{(C^1(\ooo{D}))^n} \le M$ by $H \in \HHH$ and
$M\mbox{diam}\, D < \ep M \le M$ by (3.1), the mean value theorem yields
\begin{align*}
& \max_{y\in \ooo{D}} \psi_H(y) - \min_{y\in \ooo{D}} \psi_H(y)
\le \Vert \nabla \psi_H\Vert_{(C(\ooo{D}))^n}\vert \xi-\eta\vert\\
\le& \max_{1\le k\le n} \left\Vert h_k + \sum_{j=1}^n (\ppp_kh_j)x_j
\right\Vert_{C(\ooo{D})} \max_{y,y'\in \ooo{D}} \vert y-y'\vert
\le (M + M \mbox{diam}\, D)\mbox{diam}\, D \le 2M\mbox{diam}\, D.
\end{align*}
Again condition (3.1) implies 
$$
\max_{y\in \ooo{D}} \psi_H(y) - \min_{y\in \ooo{D}} \psi_H(y)
< 2M\ep \le 2M\frac{\beta T}{4M} = \frac{1}{2}\beta T,                          \eqno{(3.2)}
$$
for all $H \in \HHH$.  By (3.2), for all $H \in \HHH$, we obtain
\begin{align*}
& \min_{y\in \ooo{D}} \va_H(y,0) - \max_{y\in \ooo{D}} \va_H(y,T)\\
= & \min_{y\in \ooo{D}} \psi_H(y) - \max_{y\in \ooo{D}} \psi_H(y) + \beta T
> -\frac{1}{2}\beta T + \beta T = \frac{1}{2}\beta T,
\end{align*}
that is,  
$$
\min_{y\in \ooo{D}} \va_H(y,0) - \max_{y\in \ooo{D}} \va_H(y,T)
> \frac{1}{2}\beta T,                       \eqno{(3.3)}
$$
for all $H \in \HHH$.

We set $\ep_0 < \frac{T}{16}$.
Then  
$$
\mbox{if $0 \le t \le 2\ep_0$, then 
$\va_H(x,t) \ge \va_H(x,2\ep_0) 
\ge \displaystyle\min_{y\in \ooo{D}} \va_H(y,0) - 2\ep_0\beta$}
$$
and
$$
\mbox{if $T-2\ep_0 \le t \le T$, then 
$\va_H(x,t) \le \va_H(x,T-2\ep_0) 
\le \displaystyle\max_{y\in \ooo{D}} \va_H(y,T) + 2\ep_0\beta$.}
$$
Therefore (3.3) and $\ep_0 < \frac{T}{16}$ yield
\begin{align*}
& \min_{x\in \ooo{D},0\le t \le 2\ep_0} \va_H(x,t)
- \max_{x\in \ooo{D},T-2\ep_0 \le t \le T} \va_H(x,t)\\
\ge & \min_{y\in \ooo{D}} \va_H(y,0) - 2\ep_0\beta
- (\max_{y\in \ooo{D}} \va_H(y,T) + 2 \ep_0\beta) \\
>& \frac{1}{2}\beta T - 4\ep_0\beta > \frac{1}{4}\beta T.
\end{align*}
Setting
$$
\sigma_1(H) = \min_{x\in \ooo{D},0\le t \le 2\ep_0} \va_H(x,t),
\quad \sigma_2(H) = \max_{x\in \ooo{D},T-2\ep_0 \le t \le T} \va_H(x,t),
$$
we have
$$
\sigma_1(H) - \sigma_2(H) > \frac{1}{4}\beta T   \eqno{(3.4)}
$$
for all $H \in \HHH$.
Now we define a cut-off function $\chi \in C^{\infty}(\R)$ satisfying
$0 \le \chi \le 1$ and
$$
\chi(t) = 
\left\{ \begin{array}{rl}
&1, \quad 0 < t < T-2\ep_0,\\
&0, \quad T-\ep_0 < t < T.
\end{array}\right.
                                  \eqno{(3.5)}
$$

{\bf Second Step: Application of Carleman estimate.}\\
In this step, we apply the Carleman estimate to complete the proof.

Setting
$$\left\{\begin{array}{rl}
& v_k = u(H_1,a_k) - u(H_2,a_k), \quad R_k = u(H_2,a_k), \quad
k=1,2,..., n,\\
& H = H_2-H_1,
\end{array}\right.
$$
we have
$$\left\{\begin{array}{rl}
& \ppp_tv_k + (H_1(x)\cdot \nabla v_k) + p(x)v_k
= (H(x)\cdot \nabla R_k(x,t)) \quad
\mbox{in $Q:= D\times (0,T)$},\\
& v_k(x,0) = 0 \quad \mbox{in $D$, $k=1,2,...,n$}.
\end{array}\right.                                   \eqno{(3.6)}
$$
Moreover putting
$$
w_k = \chi\ppp_tv_k, \quad k=1,2,..., n,
$$
we obtain
$$\left\{\begin{array}{rl}
& \ppp_tw_k + (H_1(x)\cdot \nabla w_k) + p(x)w_k 
= \chi(H(x)\cdot \nabla \ppp_tR_k(x,t))
+ \chi'(t)\ppp_tv_k  \quad \mbox{in $Q$},\\
& w_k(x,0) = (H(x)\cdot \nabla a_k) \quad \mbox{in $D$, $k=1,2,...,n$}.
\end{array}\right.                                   \eqno{(3.7)}
$$
Here and henceforth we write $\chi'(t) = \frac{d\chi}{dt}(t)$.
For calculating $w_k(x,0)$, we used (3.6) and $R_k(x,0) = 
u(H_2,a_k)(x,0) = a_k(x)$.

By Lemma 2 we have
$$
(H_1(x) \cdot \nu(x)) \le 0, \quad x \in \ppp D \setminus \Gamma.
$$
Therefore, since $w_k(\cdot,T) = 0$ by (3.5), in terms of Lemma 1 by 
(3.1), we apply Proposition 1 to (3.7)
and we obtain
$$
s\int_D \vert w_k(x,0)\vert^2e^{2s\va_{H_1}(x,0)} dx
+ s^2\int_Q \vert w_k\vert^2 e^{2s\va_{H_1}(x,t)} dxdt
$$
$$
\le C\int_Q \vert \chi(H(x)\cdot \nabla \ppp_tR_k)\vert^2
e^{2s\va_{H_1}} dxdt
+ C\int_Q \vert \chi'(t)\ppp_tv_k\vert^2 e^{2s\va_{H_1}} dxdt
+ Ce^{Cs}d^2                             \eqno{(3.8)}
$$
for all $s>s_0$.  Here and henceforth we set 
$$
d = \sum_{k=1}^n \Vert w_k\Vert_{L^2(\Gamma\times (0,T))}
\le \sum_{k=1}^n \Vert \ppp_tu(H_1,a_k) - \ppp_tu(H_2,a_k)
\Vert_{L^2(\Gamma\times (0,T))}.
$$
Since $\chi'(t) = 0$ $\forall t\in [0,T-2\ep_0]\cup[T-\ep_0,T]$,
we have
$$
\int_Q \vert \chi'(t)\ppp_tv_k\vert^2 e^{2s\va_{H_1}} dxdt
= \int^{T-\ep_0}_{T-2\ep_0} \int_D \vert \chi'(t)\vert^2
\vert \ppp_tu(H_1,a_k) - \ppp_tu(H_2,a_k)\vert^2 e^{2s\va_{H_1}} dxdt
$$
$$
\le Ce^{2s\sigma_2(H_1)} (\Vert \ppp_tu(H_1,a_k) \Vert_{L^2(Q)}^2
+ \Vert \ppp_tu(H_2,a_k)\Vert^2_{L^2(Q)}) 
\le Ce^{2s\sigma_2(H_1)} M_0^2.                 \eqno{(3.9)}
$$
Substituting the initial value in (3.7) into (3.8) and using the 
boundedness condition in $\HHH$ given by (1.3), we see
\begin{align*}
& s\int_D \vert (\nabla a_1(x), ..., \nabla a_n(x))^T H(x)\vert^2 
e^{2s\va_{H_1}(x,0)} dx\\
\le& C\int^T_0 \int_D \vert H(x)\vert^2 M_0^2 e^{2s(\psi_{H_1}(x)-\beta t)}dxdt
+ Ce^{2s\sigma_2(H_1)} M_0^2 + Ce^{Cs}d^2
\end{align*}
for all $s>s_0$.
By (1.4) we have
$$
\vert (\nabla a_1(x), ..., \nabla a_n(x))^T H(x)\vert
\ge C_0\vert H(x)\vert, \quad x \in \ooo{D}
$$
with some constant $C_0 > 0$.  Hence
\begin{align*}
& s\int_D \vert H(x)\vert^2 e^{2s\va_{H_1}(x,0)} dx\\
\le& C\int_D \vert H(x)\vert^2 e^{2s\va_{H_1}(x,0)} dx
+ Ce^{2s\sigma_2(H_1)} M_0^2 + Ce^{Cs}d^2
\end{align*}
for all $s>s_0$.
Taking $s_1>0$ sufficiently large, we absorb the first term on the 
right-hand side into the left-hand side, we obtain
$$
s\int_D \vert H(x)\vert^2 e^{2s\va_{H_1}(x,0)} dx
\le Ce^{2s\sigma_2(H_1)} M_0^2 + Ce^{Cs}d^2
$$
for all $s>s_1$.
Since $\va_{H_1}(x,0) \ge \sigma_1(H_1)$ for $x \in \ooo{D}$, we have
$$
\Vert H\Vert^2_{L^2(D)} \le Ce^{-2s(\sigma_1(H_1)-\sigma_2(H_1))} M_0^2 
+ Ce^{Cs}d^2
$$
for all $s>s_1$.
The inequality (3.4) implies 
$$
\Vert H\Vert^2_{L^2(D)} \le Ce^{-\frac{1}{2}\beta T s}M_0^2 
+ Ce^{Cs}d^2                                     \eqno{(3.10)}
$$
for all $s>s_1$.  
Replacing $C$ by $Ce^{Cs_1}$, we see that (3.10) holds for all 
$s>0$.
We then choose $s>0$ so that the right-hand side of (3.10) is smaller.
We consider two cases separately.
\\
{\bf Case 1: $M_0>d$.}\\
Choosing $s>0$ such that 
$$
e^{-\frac{1}{2}\beta T s}M_0^2 = e^{Cs}d^2,
$$
that is,
$$
s = \frac{2}{C+\frac{1}{2}\beta T}\log \frac{M_0}{d} > 0
$$ 
since $\frac{M_0}{d} > 1$.  Then (3.10) gives
$$
\Vert H\Vert^2_{L^2(D)} \le 2CM_0^{\frac{4C}{2C+\beta T}}
d^{\frac{2\beta T}{2C+\beta T}}.
$$
We set $\theta = \frac{\beta T}{2C+\beta T} \in (0,1)$.
\\
{\bf Case 2: $M_0\le d$.}\\
We readily see that $\Vert H\Vert^2_{L^2(D)} \le 2Ce^{Cs}d^2$.
\\

Hence, combining Cases 1 and 2, we reach 
$$
\Vert H\Vert_{L^2(D)} \le C(d^{\theta} + d).
$$
Thus the proof of Theorem 1 is complete.
\section{Proof of Theorem 2}

The proof is similar to Theorem 1, and simpler because the principal part,
i.e. $H(x),$ is known.  Setting
$$
y = u(p_1,a) - u(p_2,a), \quad f = p_1-p_2, \quad R=-u(p_2,a),
$$
we have
$$\left\{ \begin{array}{rl}
& \ppp_ty + (H(x)\cdot\nabla y) + p_1(x)y = R(x,t)f(x) \quad \mbox{in 
$Q := D\times (0,T)$}, \\
& y(x,0) = 0 \quad \mbox{in $D$}.
\end{array}\right.
                                                             \eqno{(4.1)}
$$
Let $\ep_0>0$ be chosen similarly to that of the proof of 
Theorem 1 in Section 3.
We consider $\chi \in C^{\infty}(\R),$ 
$0 \le \chi \le 1$ defined in (3.5).
Putting
$$
z = \chi\ppp_ty,
$$
we obtain
$$\left\{ \begin{array}{rl}
& \ppp_tz + (H(x)\cdot\nabla z) + p_1(x)z = \chi(\ppp_tR)(x,t)f(x)
+ \chi'(t)\ppp_ty \quad \mbox{in $Q$}, \\
& z(x,0) = R(x,0)f(x) \quad \mbox{in $D$}.
\end{array}\right.
$$
By Lemma 2, noting that $(H\cdot\nu) \le 0$ on $\ppp D \setminus \Gamma$
implies $\Gamma_+ \subset \Gamma$, we can apply Proposition 1 to the last system
$$
s\int_D \vert z(x,0)\vert^2 e^{2s\va_H(x,0)} dx 
+ s^2\int_Q \vert z\vert^2 e^{2s\va_H} dxdt
$$
$$
\le C\int_Q \vert \chi (\ppp_tR)f\vert^2 e^{2s\va_H(x,t)} dxdt
+ C\int_Q \vert \chi'\ppp_ty\vert^2 e^{2s\va_H(x,t)} dxdt
+ Ce^{Cs}d_1^2                     \eqno{(4.2)}
$$
for all $s > s_0,$ where we set 
$$
d_1 = \Vert z\Vert_{L^2(\Gamma \times (0,T)}
\le \Vert \ppp_tu(p_1,a) - \ppp_tu(p_2,a)\Vert_{L^2(\Gamma \times (0,T))}.
$$
Similarly to (3.9), keeping in mind the definition of $\sigma_2(H),$ contained in the proof of Theorem 1 in Section 3, we can estimate
$$
\int_Q \vert \chi'\ppp_ty\vert^2 e^{2s\va_H(x,t)} dxdt
= \int^{T-\ep_0}_{T-2\ep_0}\int_D 
\vert \chi'\ppp_ty\vert^2 e^{2s\va_H(x,t)} dxdt   \eqno{(4.3)}
$$
$$
\le Ce^{2s\sigma_2(H)} \int^{T-\ep_0}_{T-2\ep_0}\int_D \vert \ppp_ty\vert^2 
dx dt = Ce^{2s\sigma_2(H)}\Vert \ppp_ty\Vert^2
_{L^2(T-2\ep_0,T-\ep_0;L^2(D))}.              
$$
Therefore, by assumption (1.7), that is $\Vert \ppp_t u(p_2,a)\Vert_{L^2(0,T;L^{\infty}(D))}
\le M_0$, we obtain
$$
\int_Q \vert \chi'\ppp_ty\vert^2 e^{2s\va_H(x,t)} dxdt
\le Ce^{2s\sigma_2(H)}M_0^2.                     \eqno{(4.4)}
$$
Moreover,
$$
\int_Q \vert \chi(\ppp_tR)f\vert^2 e^{2s\va_H(x,t)} dxdt
\le CM_0^2\int_Q \vert f(x)\vert^2 e^{2s\va_H(x,t)} dx dt.
$$
By (4.1), we can calculate 
$$
z(x,0) = \chi(0)\ppp_ty(x,0) = R(x,0)f(x) = -a(x)f(x), \quad x\in D.
$$
Moreover, since 
$$
e^{2s\va_H(x,t)} \le e^{2s\va_H(x,0)} \quad \mbox{for $(x,t) \in Q$},
$$
the assumption (1.6) yields
$$
s\int_D \vert f(x)\vert^2 e^{2s\va_H(x,0)} dx
\le C\int_D \vert f(x)\vert^2 e^{2s\va_H(x,0)}dx
+ Ce^{2s\sigma_2(H)}M_0^2 + Ce^{Cs}d_1^2
$$
for all $s>s_0$.
Therefore similarly to (3.10), we can obtain
$$
\Vert f\Vert^2_{L^2(D)} \le Ce^{-\frac{1}{2}\beta Ts}M_0^2
+ Ce^{Cs}d_1^2
$$
for all $s>s_0$.  We can follow 
the proof of 
Theorem 1, and the proof of Theorem 2 is complete.
\section{Closing remarks}

In this section, we discuss the determination of both coefficients and 
solution to the transport equation also with other formulations of 
observation data.
For simplicity, we consider linearized problems of inverse problems 
of determining a zeroth order coefficient $p$ in (1.1).
We can similarly study non-linearized inverse problems of determining 
a vector-valued function $H(x)$, but we here omit details.

We set $y=u(p_1,a)-u(p_2,a)$.  Using the notations of the previous sections, 
we consider
$$
\left\{ \begin{array}{rl}
& \ppp_ty(x,t) + (H(x)\cdot\nabla y) + p_1(x)y(x,t) = R(x,t)f(x), \quad
x\in D, \, 0<t<T, \\
& y(x,0) = 0, \quad x \in D.
\end{array}\right.
                                                   \eqno{(5.1)}
$$
Henceforth we assume that $H\in \HHH$, $p_1\in L^{\infty}(D)$,
$R \in H^1(0,T;L^{\infty}(D))$ are fixed and $f \in L^2(D)$,
and $R(x,0) \ne 0$ on $\ooo{D}$.
Then, in terms of Lemma 2, by the same argument in First Step of the proof 
of Theorem 1, we can choose a domain $D$ such that 
$$\left\{ \begin{array}{rl}
& \mbox{$D$ depends on $H, T, \OOO$},\quad 
\ppp D = \ooo{\Gamma_+ \cup \Gamma_-}, \\
&\mbox{where $\Gamma_- := \{ x\in \ppp D;\, H(x)\cdot \nu(x) \le 0\},\quad
\Gamma_+ := \{ x\in \ppp D;\, H(x)\cdot \nu(x) \ge 0\} \subset \Gamma$},\\
& \mbox{$\Gamma_+ \cap \Gamma_-$ has no interior points},\quad 
x_0 \in \ppp D \cap \ppp\OOO \subset \Gamma.
\end{array}\right.
                                            \eqno{(5.2)}
$$
We note that the domain $D$ can depend on $H$ and fixing $H$, we
can consider that $D$ is fixed for our discussions in this section.
\\

{\bf \S 5.1. Determination of solution $y(x,t)$.}

In Theorems 1 and 2, we are restricted to the determination of 
coefficients, and we do not consider the determination of 
solution $u$ itself to (1.1).
For the linearized inverse problem, the same argument as the proof of 
Theorem 2 yields conditional H\"older stability in determing $f$ in $D$
and $y$ in $D\times (0,\ep)$ where $\ep>0$ is sufficiently small.
\\
{\bf Proposition 2.}\\
There exist constants $C>0$, $\theta_1 \in (0,1)$ and small 
$\ep>0$ such that 
$$
\Vert f\Vert_{L^2(D)} + \Vert y\Vert_{H^1(0,\ep;L^2(D))}
\le C(\Vert \ppp_ty\Vert^{\theta_1}_{L^2(\Gamma \times (0,T))}
+ \Vert \ppp_ty\Vert_{L^2(\Gamma \times (0,T))})
                                                   \eqno{(5.3)}
$$
provided that $\Vert \ppp_ty\Vert_{L^2(0,T;L^{\infty}(D))}\le M_0$ with 
arbitrarily fixed constant $M_0>0$.
The constants $C>0$ and $\theta_1 \in (0,1)$, $\ep>0$ depend on 
$\OOO, T, \Gamma, H$.
\\

The interval length $\ep$ for estimating $y$ is small, which means that the 
estimate is only near $t=0$.  As Example in Section 1 shows, 
we cannot prove $\ep = T$.
\\
\vspace{0.1cm}
{\bf Proof of Proposition 2.}
The estimation of $\Vert f\Vert_{L^2(D)}$ is the same as Theorem 2.
We recall that $Q=D\times (0,T)$.
Now we estimate $y$.  By the same manner as in (4.2) and (4.4), we obtain
\begin{align*}
& s^2\int_Q \vert \chi\ppp_ty\vert^2 e^{2s\va_H(x,t)} dxdt\\
\le& C\int_Q \vert \chi(\ppp_tR)f\vert^2 e^{2s\va_H(x,t)} dxdt
+ Ce^{2s\sigma_2(H)}M_0^2 + Ce^{Cs}d_1^2\\
\le &Ce^{2Cs}\Vert f\Vert^2_{L^2(D)} + Ce^{2s\sigma_2(H)}M_0^2
+ Ce^{Cs}d_1^2.
\end{align*}
Here we used $\ppp_tR \in L^2(0,T;L^{\infty}(D))$ and we set
$$
d_1 = \Vert \ppp_ty\Vert_{L^2(\Gamma \times (0,T))}.   \eqno{(5.4)}
$$
Since $\va_H(x,t) \ge \sigma_1(H)$ for $x \in D$ and $0\le t \le 2\ep_0
=: \ep$ and $\chi(t) = 1$ for $0 \le t \le 2\ep_0$, we have
$$
s^2e^{2s\sigma_1(H)}\Vert \ppp_ty\Vert^2_{L^2(0,\ep;L^2(D))}
\le Ce^{2Cs}\Vert f\Vert_{L^2(D)}^2 
+ Ce^{2s\sigma_2(H)}M_0^2 + Ce^{Cs}d_1^2.
$$
Since we have already proved $\Vert f\Vert_{L^2(D)}
\le C(d_1 + d_1^{\theta})$, we reach 
$$
\Vert \ppp_ty\Vert^2_{L^2(0,\ep;L^2(D))}
\le Ce^{2Cs}(d_1^2+d_1^{2\theta})
+ Ce^{-2s(\sigma_1(H)-\sigma_2(H))}M_0^2
$$
for all large $s>0$.  By (3.4), we can argue in the same was as
for (3.10), and we complete the estimate of $\Vert \ppp_ty\Vert
_{L^2(0,\ep;L^2(D))}$ with some $\theta_1 \in (0,\theta)$.
Since $y(\cdot,0) = 0$ in $D$, we easily prove
$\Vert y\Vert_{H^1(0,\ep;L^2(D))} \le C\Vert \ppp_ty\Vert
_{L^2(0,\ep;L^2(D))}$.  Thus the proof of Proposition 2 is complete.
\\

{\bf \S 5.2. Other formulations of the inverse problem}

So far we formulate the inverse problem with data of $y$ on 
$D \times \{0\}$ and $\Gamma_+ \times (0,T)$.
As possible data sets, we can introduce:
\begin{itemize}
\item
{\bf Case I:} $D \times \{T\}$ and $\Gamma_- \times (0,T)$
\item
{\bf Case II:} $D \times \{0\}$ and $\Gamma_- \times (0,T)$
\item
{\bf Case III:} $D \times \{T\}$ and $\Gamma_+ \times (0,T)$
\item
{\bf Case IV:} $D \times \{0,T\}$ and $\Gamma_+ \times (0,T)$
\item
{\bf Case V:} $D \times \{0,T\}$ and $\Gamma_- \times (0,T)$
\item
{\bf Case VI:} $D \times \{0\}$ and $(\Gamma_- \cup \Gamma_+) \times (0,T)$
\item
{\bf Case VII:} $D \times \{T\}$ and $(\Gamma_- \cup \Gamma_+) \times (0,T)$
\end{itemize}
We note that we exclude the two cases of data on 
$(\Gamma_- \cup \Gamma_+) \times (0,T)$ and $D \times \{0,T\}$, 
because we can expect that such two data sets are both too poor 
and that we cannot 
prove any uniqueness for the corresponding invers problems.

By the change of variables $t \longrightarrow T-t$, setting
$v(x,t) = u(x,T-t)$ for $x \in D$ and $0 < t < T$, we have
$\ppp_tv(x,t) = -\ppp_tu(x,T-t)$, so that (5.1) can be written 
in terms of $v$ and $-H(x)$.  Therefore we can exchange $\Gamma_+$ by 
$\Gamma_-$, $\Gamma_-$ by $\Gamma_+$, $u(x,0)$ by $v(x,T)$.
Hence Case I is equivalent to our currently considered case of data
on $D \times \{0\}$ and $\Gamma_+ \times (0,T)$,
Case II to Case III, Case IV to Case V, Case VI to Case VII.
Thus it is sufficient to discuss the three Cases II, IV. VI.
\\

{\bf Case II.}

In Case II, we do not know positive results for the inverse problem by 
our method.   In fact, we
cannot prove an adequate Carleman estimate which corresponds to 
Propostion 1, and for the proof of Carleman estimate in Case II,
the signs of $B_H(x,0)$ and $B_H(x)(H(x)\cdot\nu)$ on $\Gamma_- 
\times (0,T)$ are contradictive.

Data on $\Gamma_+ \times (0,T)$ may be meaningful for determining $f$.  
On the other hand, data on $\Gamma_- \times (0,T)$ are not 
meaningful for $f$, but with given $f$, data on $\Gamma_- \times (0,T)$
and $y(\cdot,0)$ perfectly determine $y$ on $D\times (0,T)$ by a usual energy 
estimate (Lemma 3 below).   
\\

{\bf Case IV.}
     
The result is different and we can 
prove an unconditional Lipschitz stability.
\\
{\bf Proposition 3.}\\
There exists a constant $C>0$ such that 
$$
\Vert f\Vert_{L^2(D)} + \Vert y\Vert_{H^1(0,T;L^2(D))}
\le C(\Vert \ppp_ty\Vert_{L^2(\Gamma \times (0,T))}
+ \Vert y(\cdot,T)\Vert_{H^1(D)}).
$$
Here the constant $C>0$ is independent of any bounds of 
$y$ and dependent on $\OOO, T, \Gamma, H$, which means that the 
stability is unconditional.
\\
{\bf Proof.}\\
Since we are given $y(\cdot,T)$, we need not the cut-off function 
$\chi$ satisfying (3.5).  Setting $y_1 = \ppp_ty$, we differentiate the
first equation in (5.1) to have
$$
\left\{ \begin{array}{rl}
& \ppp_ty_1(x,t) + (H(x)\cdot\nabla y_1) + p_1(x)y_1(x,t) 
= (\ppp_tR(x,t))f(x), \quad
x\in D, \, 0<t<T, \\
& y_1(x,0) = R(x,0)f(x), \quad x \in D
\end{array}\right.
                                                   \eqno{(5.5)}
$$
and
$$
y_1(x,T) = -(H(x)\cdot \nabla y(x,T)) - p_1(x)y(x,T) + R(x,T)f(x),
\quad x \in D.
$$
Applying Proposition 1, we obtain
\begin{align*}
&s\int_D \vert R(x,0)\vert^2\vert f(x)\vert^2 e^{2s\va_H(x,0)} dx
+ s^2\int_Q \vert \ppp_ty\vert^2 e^{2s\va_H(x,t)} dxdt\\
\le& C\int_Q \vert \ppp_tR\vert^2 \vert f\vert^2e^{2s\va_H(x,t)} dxdt
+ Cs\int_D \vert (H(x)\cdot \nabla y(x,T)) 
+ p_1(x)y(x,T)\vert^2 e^{2s\va_H(x,T)} dx\\
+& Cs\int_D \vert R(x,T)\vert^2 \vert f\vert^2 e^{2s\va_H(x,T)} dx
+ Ce^{Cs}d_1^2
\end{align*}
for all large $s>0$.
Here we recall that $d_1$ is defined by (5.4).

By $\ppp_tR \in L^2(0,T;L^{\infty}(D))$ and $R(x,0) \ne 0$ for 
$x \in \ooo{D}$ , since $\va_H(x,t) \le \va_H(x,0)$ for $(x,t) \in Q$
and $\sigma_2(H) \ge\displaystyle \max_{x\in \ooo{D}} \va_H(x,T),$ we obtain
\begin{align*}
&s\int_D \vert f(x)\vert^2 e^{2s\va_H(x,0)} dx
+ s^2e^{-2sC_1}\Vert \ppp_ty\Vert^2_{L^2(0,T;L^2(D))}\\
\le& C\int_D \vert f(x)\vert^2e^{2s\va_H(x,0)} dx
+ Cse^{2sC_1}\Vert y(\cdot,T)\Vert^2_{H^1(D)}
+ Cse^{2s\sigma_2(H)}\Vert f\Vert_{L^2(D)}^2 + Ce^{Cs}d_1^2.
\end{align*}
Here we set $C_1 = \Vert\va_H\Vert_{L^{\infty}(Q)}$.
Choosing $s>0$ large, we can absorb the first term on the right-hand side
into the left-hand side.  Using
$$
\va_H(x,0) \ge \min_{x\in \ooo{D}} \va_H(x,0) \ge \sigma_1(H),
$$
and dividing by $e^{2s\sigma_1(H)}$, we reach
\begin{align*}
& \Vert f\Vert^2_{L^2(D)}
+ s^2e^{-2s(C_1+\sigma_1(H))}\Vert \ppp_ty\Vert^2_{L^2(0,T;L^2(D))}\\
\le& Cse^{2sC_1}\Vert y(\cdot,T)\Vert^2_{H^1(D)}
+ Cse^{-2s(\sigma_1(H)-\sigma_2(H))}\Vert f\Vert_{L^2(D)}^2 + Ce^{Cs}d_1^2
\end{align*}
for all large $s>0$.
In terms of (3.4), again choosing $s>0$ large, we can absorb the second term 
on the right-hand side into the left-hand side, and we complete the proof of
Proposition 3.
\\
  
{\bf Case VI.}

Also in this case, we can prove an unconditional Lipschitz stability 
estimate but the norm of $y$ is not the same.
\\
{\bf Proposition 4.}\\
There exists a constant $C>0$ such that 
$$
\Vert f\Vert_{L^2(D)} + \Vert y\Vert_{W^{1,\infty}(0,T;L^2(D))}
\le C\Vert \ppp_ty\Vert_{L^2(\ppp D \times (0,T))}
$$
with some constant $C>0$ which is independent of any bounds of 
$\Vert \ppp_ty\Vert_{L^2(0,T;L^{\infty}(D))}$.
\\
{\bf Proof.}\\
With the same notation $z = \chi\ppp_ty$, in the proof of Theorem 2,
we can obtain (4.2) and (4.3), where $d_1$ is defined by 
$$
d_1 = \Vert \chi\ppp_ty\Vert_{L^2(\Gamma_+\times (0,T))}.
                                                   \eqno{(5.5)}
$$
Therefore
\begin{align*}
&s\int_D \vert \ppp_ty(x,0)\vert^2 e^{2s\va_H(x,0)} dx
\le C\int_Q \vert f(x)\vert^2 e^{2s\va_H(x,t)} dxdt\\
+ & Ce^{2s\sigma_2(H)}\Vert \ppp_ty\Vert^2
_{L^2(T-2\ep_0,T-\ep_0;L^2(D))}
+ Ce^{Cs}d_1^2
\end{align*}
for all large $s>0$.
Since $\ppp_ty(x,0) = R(x,0)f(x)$ and 
$\va_H(x,t) \le \va_H(x,0)$ for $(x,t) \in Q$, we have
\begin{align*}
&s\int_D \vert f(x)\vert^2 e^{2s\va_H(x,0)} dx
\le C\int_D \vert f(x)\vert^2 e^{2s\va_H(x,0)} dxdt\\
+ & Ce^{2s\sigma_2(H)}\Vert \ppp_ty\Vert^2
_{L^2(T-2\ep_0,T-\ep_0;L^2(D))}
+ Ce^{Cs}d_1^2
\end{align*}
Choosing $s>0$ large, we can absorb the first term on the right-hand side 
into the left-hand side, we obtain
$$
\int_D \vert f(x)\vert^2 e^{2s\va_H(x,0)} dx
\le Ce^{2s\sigma_2(H)}\Vert \ppp_ty\Vert^2_{L^2(0,T;L^2(D))}
+ Ce^{Cs}d_1^2                            \eqno{(5.6)}
$$

The following lemma is usual energy estimation and estimates
$\Vert \ppp_ty\Vert^2_{L^2(0,T;L^2(D))}$.
\\
{\bf Lemma 3.}\\
There exists a constant $C>0$ depending on $R$, such that  
$$
\Vert \ppp_ty(\cdot,t)\Vert_{L^2(D)}
\le C(\Vert f\Vert_{L^2(D)}
+ \Vert \ppp_ty\Vert_{L^2(\Gamma_- \times (0,T))}), \quad 
0\le t\le T.
$$
\\
{\bf Proof of Lemma 3.}\\
Setting $y_1 = \ppp_ty$, we have (5.5).
Multiplying the first equation in (5.5) with $2y_1$ and integrating 
in $D$, we have
\begin{align*}
& \ppp_t\int_D \vert y_1(x,t)\vert^2 dx 
+ 2\int_D \sum_{j=1}^n h_j(x)(\ppp_jy_1(x,t))y_1(x,t) dx 
+ 2\int_D p_1(x)\vert y_1(x,t)\vert^2 dx\\
=& 2\int_D (\ppp_tR)(x,t)f(x)y_1(x,t) dx.
\end{align*}
Here the integration by parts provides
\begin{align*}
&2\int_D \sum_{j=1}^n h_j(x)(\ppp_jy_1(x,t))y_1(x,t) dx
= \int_D \sum_{j=1}^n h_j(x)\ppp_j(\vert y_1(x,t)\vert^2) dx\\
=& \int_{\ppp D} (H\cdot\nu) \vert y_1\vert^2 dS
- \int_D (\mbox{div}\, H)\vert y_1(x,t)\vert^2 dx\\
\ge &\int_{\Gamma_-} (H\cdot\nu) \vert y_1\vert^2 dS
- \int_D (\mbox{div}\, H)\vert y_1(x,t)\vert^2 dx.
\end{align*}
Setting $E(t) = \int_D \vert y_1(x,t)\vert^2 dx$, the Cauchy-Schwarz 
inequality yields
\begin{align*}
& \frac{dE(t)}{dt} - C\int_{\Gamma_-} \vert y_1\vert^2 dS
\le CE(t) + \int_D \vert (\ppp_tR)f\vert^2 dx + \int_D \vert y_1\vert^2 dx\\
\le& CE(t) + C\Vert f\Vert^2_{L^2(D)}.
\end{align*}
Therefore by integration, we have
$$
E(t) \le E(0) + C\int^t_0 E(\eta) d\eta 
+  C(\Vert y_1\Vert^2_{L^2(\Gamma_-\times (0,T))} 
+ \Vert f\Vert^2_{L^2(D)}), \quad 0\le t\le T.
$$
Since $E(0) = \int_D \vert R(x,0)f(x)\vert^2 dx \le C\Vert f\Vert^2_{L^2(D)}$,
the Gronwall inequality completes the proof of Lemma 3.
\\

Now we return to the proof of Proposition 4.
Applying Lemma 3 in the first term on the right-hand side of (5.6),
we have
\begin{align*}
& \int_D \vert f(x)\vert^2 e^{2s\va_H(x,0)} dx
\le Ce^{2s\sigma_2(H)}\Vert f\Vert^2_{L^2(D)}\\
+ &Ce^{2s\sigma_2(H)}\Vert \ppp_ty\Vert^2_{L^2(\Gamma_- \times (0,T))}
+ Ce^{Cs}\Vert \ppp_ty\Vert^2_{L^2(\Gamma_+ \times (0,T))}
\end{align*}
for all large $s>0$.
Noting (3.4) and using $\va_H(x,0) \ge \sigma_1(H)$ for $x \in D$,
we divide both sides with $e^{2s\sigma_1(H)}$, we obtain
$$
\Vert f\Vert^2_{L^2(D)}
\le Ce^{-\frac{1}{2}\beta Ts}\Vert f\Vert^2_{L^2(D)}
+ Ce^{Cs}\Vert \ppp_ty\Vert^2_{L^2((\Gamma_+ \cup \Gamma_-)\times 
(0,T))}
$$
for all large $s>0$.  Again choosing $s>0$ sufficiently large, we can 
absorb the first term on the right-hand side into the left-hand side, so 
that 
$$
\Vert f\Vert^2_{L^2(D)} \le Ce^{Cs}\Vert \ppp_ty\Vert^2
_{L^2(\ppp D\times (0,T))}.
$$
With Lemma 3, we see also 
$$
\Vert \ppp_ty\Vert_{L^{\infty}(0,T;L^2(D))}
\le C\Vert \ppp_ty\Vert_{L^2(\ppp D\times (0,T))}.
$$
Thus the proof of Proposition 4 is completed.\\

{\bf Acknowledgement.} 
This work was partially supported by Grant-in-Aid for Scientific Research (S) 15H05740 of
Japan Society for the Promotion of Science,
and by NSFC (No. 11771270, 91730303) and the 
\lq\lq RUDN University Program 5-100''.  
The first and second authors were visitors at 
The University of Tokyo in February 2018, supported by the above grant of the Japan Society for the Promotion of Science.
Most of the paper was completed during 
the stay of the fourth author as GNAMPA visitor at the IN$\delta$AM Unit of the Universit\`a degli Studi di Napoli Federico 
II in January 2019 and he is grateful for the supports. 
The fourth author was a Visiting Scholar at Rome in April 2018 supported by the University of Rome \lq\lq 
Tor Vergata''. This work was supported also by the Istituto Nazionale di Alta Matematica 
(INdAM), through the GNAMPA Research Project 2017 \lq\lq 
Comportamento asintotico e controllo di equazioni di evoluzione non lineari''. 
Moreover, this research was 
performed within the framework of the
French-German-Italian Laboratoire International Associ\'e (LIA), named COPDESC, on Applied Analysis, issued by CNRS, MPI and INdAM. This work was also supported by 
the research project of the Universit\`a di Napoli Federico II: \lq\lq Spectral 
and Geometrical Inequalities''.

\end{document}